\documentclass[12pt]{article}
\usepackage[T2A]{fontenc}
\usepackage[utf8]{inputenc}
\usepackage{amssymb,amsmath,amsfonts,amsthm,amscd,latexsym,indentfirst}
\def\id{\textrm{id}}

\begin{document}

\hfill{11B68, 16W99 (MSC2010)}

\begin{center}
{\Large
Rota---Baxter operators and Bernoulli polynomials}

Vsevolod Gubarev
\end{center}

\begin{abstract}
We develop the connection between Rota---Baxter operators 
arisen from algebra and mathematical physics and Bernoulli polynomials.
We state that a trivial property of Rota---Baxter operators implies
the symmetry of the power sum polynomials and Bernoulli polynomials.
We show how Rota---Baxter operators equalities rewritten in terms of 
Bernoulli polynomials generate identities for the last.

\end{abstract}

\section{Introduction}

Given an algebra $A$ and a scalar $\lambda\in F$, where $F$ is a~ground field, 
a~linear operator $R\colon A\rightarrow A$ is called a Rota---Baxter operator 
on $A$ of weight~$\lambda$ if the following identity
\begin{equation}\label{RB}
R(x)R(y) = R( R(x)y + xR(y) + \lambda xy )
\end{equation}
holds for all $x,y\in A$. The algebra $A$ is called Rota---Baxter algebra.
By algebra we mean a vector space endowed a bilinear not necessary associative product. 

The notion of Rota---Baxter operator was introduced 
by G. Baxter~\cite{Baxter} in 1960 as formal generalization 
of integration by parts formula (when $\lambda = 0$) and then developed 
by G.-C. Rota~\cite{Rota} and others \cite{Atkinson,Cartier}.

In 1980s, the deep connection between constant solutions 
of the classical Yang---Baxter equation from mathematical 
physics and Rota---Baxter operators of weight zero on a semi\-simple 
finite-dimensional Lie algebra was discovered \cite{BelaDrin82,Semenov83}.

To the moment, applications of Rota---Baxter operators 
in symmetric polynomials, quantum field renormalization, 
Loday algebras, shuffle algebra etc were found \cite{Atkinson,Renorm,GuoMonograph,Aguiar00,Fard02,GubKol2013,Shuffle}.
The notion of Rota---Baxter operator is useful in such branch of number theory as multiple zeta function~\cite{MZV,Zhao}.

In 1966, J. Miller found an interesting connection between 
Rota---Baxter operators and the power sum polynomials~\cite{Miller}.
We start with an algebra $A$ which is unital and power-associative 
(it means every one-generated subalgebra is associative).
Let $R$ be a Rota---Baxter operator on $A$ of weight~$-1$. 
Denote by~1 the unit of $A$ and put $a = R(1)$.
Then $R(a^n) = F_n(a)$ for all $n\in\mathbb{N}$, where $F_n(m) = \sum\limits_{j=1}^{m}j^n$.

In 2010, O. Ogievetsky and V. Schechtman restated this connection 
to find a new proof of the Schl\"{o}milch---Ramanujan formula~\cite{Ogievetsky}.
In 2017, the author reproved the connection formula to apply for Rota---Baxter operators 
of nonzero weight on the matrix algebra~\cite{Unital}.

Our goal is to develop this connection. 
There exist several different proofs of the symmetry of the power sum polynomials 
$F_n(y) = (-1)^{n+1}F_n(-1-y)$ and the symmetry of Bernoulli polynomials 
$B_n(x) = (-1)^n B_n(1-x)$ involving infinite series, generating functions 
or some special identities~\cite{Lehmer,Newsome,Tuenter}. 
In~\S2, we prove that the trivial property that
the linear operator $(\id-R)$ is a Rota---Baxter operator (provided that $R$ is so)
implies both symmetries.

In~\S3, we show how identities concerned Rota---Baxter operators 
rewritten in terms of Bernoulli polynomials and Bernoulli numbers 
generate a plenty of identities for both of them.
In particular, we find a symmetric expression for the product $B_i(x)B_j(x)B_k(x)$
and count the sum $\sum\limits_{i+j+k=n;i,j,k>0}\mathcal{B}_i(x)\mathcal{B}_j(x)\mathcal{B}_k(x)$,
where $\mathcal{B}_s(x) = B_s(x)/s$ is the divided Bernoulli polynomial. The approach for counting the same sum for usual
(not divided) Bernoulli polynomials was developed in~\cite{KimProducts}.
About the products of Bernoulli polynomials and numbers see also~\cite{Carlitz,Gessel,AgohDilcher2}.

\section{Symmetry of the power sum polynomials}

{\bf Statement 1}~\cite{GuoMonograph}.
Given an RB-operator $P$ of weight $\lambda$,

a) the operator $-P-\lambda\id$ is an RB-operator of weight $\lambda$,

b) the operator $\lambda^{-1}P$ is an RB-operator of weight 1, provided $\lambda\neq0$.

Given an algebra $A$, let us define a map $\phi$ on the set of all RB-operators on $A$
as $\phi(P)=-P-\lambda(P)\id$. It is clear that $\phi^2$ coincides with the identity map.

Let $F_n(m) = \sum\limits_{j=1}^{m}j^n$ for natural $n,m$.
Bernoulli polynomials $B_n(x)$ are connected with the power sum polynomials in the following way:
\begin{equation}\label{SumsViaBernoulliPolyn}
F_n(m) = \frac{B_{n+1}(m+1)-B_{n+1}}{n+1}.
\end{equation}

{\bf Statement~2}~\cite{Miller,Ogievetsky,Unital}.
Let $A$ be a unital power-associative algebra, $R$ be an RB-operator 
on $A$ of weight~$\lambda$, $a = R(1)$. 
Then $R(a^n) = (-\lambda)^{n+1}F_n(-a/\lambda)$ for all $n\in\mathbb{N}$.
In particular, $R(a^n) = F_n(a)$ for all $n\in\mathbb{N}$ provided that $\lambda = -1$.

Let us show how the trivial property of Rota---Baxter operators from Statement~1a
implies the symmetry of the power sum polynomials and the symmetry of Bernoulli polynomials.

{\bf Lemma}. 
Let $A$ be a unital power-associative algebra, 
$R$ be an RB-operator on $A$ of weight~$-1$, $a = R(1)$ and $b = \phi(R)(1) = 1-a$.
For all positive natural $n$, we have
\begin{equation}\label{InterestEqRB}
R(a^n) - a^n = (-1)^{n+1}(\phi(R)(b^n)-b^n).
\end{equation}

\begin{proof}
From
$$
(-1)^{n+1}(\phi(R)(b^n)-b^n)
 = (-1)^{n+1}(-R)(b^n)
 = R((-b)^n) 
 = R((a-1)^n),
$$
we conclude that it is enough to state $R((a-1)^n) = R(a^n)-a^n$.
We prove the last equality by induction on~$n$. For $n = 1$,
we get the true equality $\frac{a^2-a}{2} = \frac{a^2-a}{2}$.
Suppose that we have done for all natural numbers less than~$n$.
Now we rewrite $R((a-1)^{n+1})$ by~\eqref{RB} and the induction hypothesis
\begin{multline}\label{InductStepMain}
R((a-1)^{n+1})
 = R((a-1)^n(a-1))
 = R((a-1)^nR(1)) - R((a-1)^n) \\
 = R((a-1)^n)R(1) - R(R((a-1)^n)) + R((a-1)^n) - R((a-1)^n) \\
 = (R(a^n)-a^n)a - R(R(a^n)-a^n).
\end{multline}
Again by~\eqref{RB}, we calculate
\begin{multline}\label{InductStepAdd}
R(R(a^n))
 = R(R(a^n)\cdot 1)
 = R(a^n)R(1) - R(a^nR(1)) + R(a^n\cdot1) \\
 = R(a^n)a - R(a^{n+1}) + R(a^n).
\end{multline}
Substituting~\eqref{InductStepAdd} in~\eqref{InductStepMain}
gives us the proof of the inductive step.
\end{proof}

{\bf Theorem 1}. Let $n$ be a positive natural number. Then 

a) $F_n(y) = (-1)^{n+1}F_n(-1-y)$ for all $y$,

b) $B_n(x) = (-1)^n B_n(1-x)$ for all $x$.

\begin{proof}
a) Let us consider a unital power-associative algebra~$A$
with a Rota---Baxter operator~$R$ on~$A$ of weight~$\lambda = -1$.
Put $a = R(1)$. Define $Q = \phi(R) = \id - R$ and $b = Q(1) = 1-a$.
Applying Statement~2 to the formula~\eqref{InterestEqRB},
we get $F_n(a-1) = (-1)^{n+1}F_n(-a)$, which gives a).

b) It follows from a) via~\eqref{SumsViaBernoulliPolyn}.
\end{proof}

\section{Product of two Bernoulli polynomials}

For any $n$,
\begin{equation}\label{sum-formula}
F_n(m) = \frac{1}{n+1}\sum\limits_{j=0}^n(-1)^j 
 \binom{n+1}{j} B_j m^{n+1-j},
\end{equation}
where $B_j$ is Bernoulli number.

Let us show how Rota---Baxter operators could generate a plenty of identities for Bernoulli numbers and Bernoulli polynomials.
Let $A$ be a power-associative algebra and $R$ be a Rota---Baxter operator of weight~$-1$ on $A$, $a = R(1)$.
Consider the equality
\begin{equation}\label{RBnm}
R(a^n)R(a^m) = R(R(a^n)a^m + a^nR(a^m) - a^{n+m}),\quad n,m\in \mathbb{N}.
\end{equation}
The LHS of~\eqref{RBnm} by~\eqref{SumsViaBernoulliPolyn} and Statement~2 is equal to
\begin{multline}\label{RBnmLeft}
R(a^n)R(a^m) 
 = \mathcal{B}_{n+1}(a+1)\mathcal{B}_{m+1}(a+1)
 -\mathcal{B}_{m+1}\mathcal{B}_{n+1}(a+1)
 -\mathcal{B}_{n+1}\mathcal{B}_{m+1}(a+1) \\
 +\mathcal{B}_{n+1}\mathcal{B}_{m+1},
\end{multline}
where $\mathcal{B}_n(x) = B_n(x)/n$ and $\mathcal{B}_n = B_n/n$.

Let us write down the RHS of~\eqref{RBnm} by~\eqref{SumsViaBernoulliPolyn},~\eqref{sum-formula} and Statement~2,
\begin{multline}\label{RBnmRight}
R(R(a^n)a^m + a^nR(a^m) - a^{n+m}) \\
 = \sum\limits_{i=0}^n (-1)^{n-i}\frac{1}{n+1}\binom{n+1}{n-i}B_{n-i}(\mathcal{B}_{m+2+i}(a+1)-\mathcal{B}_{m+2+i} ) \\
 + \sum\limits_{j=0}^m (-1)^{m-j}\frac{1}{m+1}\binom{m+1}{m-j}B_{m-j}(\mathcal{B}_{n+2+j}(a+1)-\mathcal{B}_{n+2+j} ) \\
 -\mathcal{B}_{n+m+1}(a+1)+\mathcal{B}_{n+m+1}.
\end{multline}
Comparing~\eqref{RBnmLeft} and~\eqref{RBnmRight}, we get the identity
\begin{equation}\label{almostNielsen}
\mathcal{B}_{i}(x)\mathcal{B}_{j}(x) - \mathcal{B}_{i}\mathcal{B}_{j}
 = \!\!\sum\limits_{l\geq0}\left(\frac{1}{i}\binom{i}{2l}+\frac{1}{j}\binom{j}{2l}\right)
 B_{2l}(\mathcal{B}_{i+j-2l}(x)-\mathcal{B}_{i+j-2l}).
\end{equation}
Here $i = n+1\geq1$, $j = m+1\geq1$ and $x = a+1$.

Up to constant, the equality~\eqref{almostNielsen} coincides with the famous identity 
\begin{equation}\label{Nielsen}
\mathcal{B}_{i}(x)\mathcal{B}_{j}(x)
 = \!\!\sum\limits_{l\geq0}\left(\frac{1}{i}\binom{i}{2l}+\frac{1}{j}\binom{j}{2l}\right)
 B_{2l}\mathcal{B}_{i+j-2l}(x) 
 + \frac{(-1)^{i-1} (i-1)!(j-1)!}{(i+j)!}B_{i+j}
\end{equation}
known at least since 1923~\cite{Nielsen}. 

{\bf Remark 1}. Writting down~\eqref{RBnm} on the first power of $a$, we get the identity
\begin{equation}\label{Agoh}
B_{n+m} 
 + \frac{1}{n+1}\sum\limits_{k=0}^n\binom{n+1}{n-k}B_{n-k}B_{m+k+1}
 + \frac{1}{m+1}\sum\limits_{l=0}^m\binom{m+1}{m-l}B_{m-l}B_{n+l+1} = 0
\end{equation}
discovered by T. Agoh in 1988~\cite{Agoh88}.

{\bf Remark 2}. 
Let us sum~\eqref{Nielsen} for $i+j = N\geq2$ and $i=1,\ldots,N-1$:
\begin{multline}\label{GesselIdReproof}
\sum\limits_{i+j=N;i,j>0}\mathcal{B}_i(x)\mathcal{B}_j(x) \\
 = \sum\limits_{i=1}^{N-1}\left(\frac{1}{i}+\frac{1}{N-i}\right)\mathcal{B}_N(x) 
 + \sum\limits_{\substack{i+j=N \\ i,j>0}}\sum\limits_{l>0}\left(\frac{1}{i}\binom{i}{2l}+\frac{1}{j}\binom{j}{2l}\right)
 B_{2l}\mathcal{B}_{N-2l}(x) \\
 + \frac{B_N}{N(N-1)}\sum\limits_{i=1}^{N-1}\frac{(-1)^{i-1} (i-1)!(N-1-i)!}{(N-2)!} \\
 = 2H_{N-1}\mathcal{B}_N(x)
 + 2\sum\limits_{l=1}^{[\frac{N-1}{2}]}B_{2l}\mathcal{B}_{N-2l}(x)\frac{1}{2l}\sum\limits_{i=1}^{N-1}\binom{i-1}{2l-1}
 + \frac{B_N}{N(N-1)}\sum\limits_{p=0}^{N-2}\frac{(-1)^p}{\binom{N-2}{p}} \\
 = 2H_{N-1}\mathcal{B}_N(x) + \frac{2}{N}\sum\limits_{l=1}^{[\frac{N-1}{2}]}\binom{N}{2l}\mathcal{B}_{2l}B_{N-2l}(x)
 + \frac{2\mathcal{B}_N}{N} \\
 = 2H_{N-1}\mathcal{B}_N(x) + \frac{2}{N}\sum\limits_{k=1}^N\binom{N}{k}\mathcal{B}_k B_{N-k}(x) + B_{N-1}(x),
\end{multline}
where $H_i = 1 + 1/2 + \ldots + 1/i$. We have used the equality~(14) from~\cite{Sury}
\begin{equation}\label{AlternOverBinom}
\sum\limits_{r=0}^n(-1)^r/\binom{n}{r} = (1+(-1)^n)\frac{n+1}{n+2}.
\end{equation}

Thus, we got in~\eqref{GesselIdReproof} the identity found by I.~Gessel in 2005~\cite{Gessel} (see also \cite{Zagier}),
\begin{equation}\label{GesselId}
\frac{N}{2}\left(-B_{N-1}(x) + \sum\limits_{k=1}^{N-1}\mathcal{B}_k(x)\mathcal{B}_{N-k}(x)\right)
 = \sum\limits_{k=1}^N\binom{N}{k}\mathcal{B}_k B_{N-k}(x) + H_{N-1}B_N(x).
\end{equation}

By the same strategy, we can compute
\begin{equation}\label{Agoh14}
\sum\limits_{k=0}^N B_k(x)B_{N-k}(x) \\
 = \frac{2}{N+2}\sum\limits_{t\geq0}\binom{N+2}{2t+2}B_{2t}B_{N-2t}(x),
\end{equation}
the identity obtained by D. Kim et al. in 2012~\cite{KimTwo} (see also~\cite{Agoh14}). 

The case $x = 0$ for~\eqref{GesselId} implies the famous identity of H. Miki~\cite{Miki} (1978) 
\begin{equation}\label{Miki}
\sum\limits_{k=2}^{N-2}\mathcal{B}_k\mathcal{B}_{N-k}
 = \sum\limits_{k=2}^{N-2}\binom{N}{k}\mathcal{B}_k \mathcal{B}_{N-k} + 2H_N \mathcal{B}_N
\end{equation}
and for~\eqref{Agoh14} it implies the identity of Yu. Matiyasevich~\cite{Matiyasevich} (1997)
\begin{equation}\label{Matiyasevich}
(N+2)\sum\limits_{k=2}^{N-2} B_k B_{N-k}
 = 2\sum\limits_{k=2}^{N-2}\binom{N+2}{k}B_k B_{N-k} + N(N+1)B_N.
\end{equation}

\section{Product of three Bernoulli polynomials}

We also may produce another identities involving the products of three, four etc Bernoulli numbers.
To do this, it is enough to consider the equality
\begin{multline}\label{RBnml}
R(a^n)R(a^m)R(a^l) 
 = R( R(a^n)R(a^l)a^m + R(a^m)R(a^l)a^n + R(a^n)R(a^m)a^l \\
 - R(a^n)a^{m+l} - R(a^m)a^{n+l} - R(a^l)a^{n+m} +a^{n+m+l})
\end{multline}
and the same equalities for four, five etc multipliers (see the formulas in~\cite{GubKol2013}).

Let us derive the explicit identity which follows from~\eqref{RBnml}.

{\bf Theorem 2}.
The following identity holds for all $i,j,k>0$,
\begin{multline}\label{SingleTripleBern}
\mathcal{B}_i(x)\mathcal{B}_j(x)\mathcal{B}_k(x) 
 = \sum\limits_{q,t\geq0}B_{2q}B_{2t-2q} \left[   
   \binom{i+j-2q}{2t-2q}\frac{1}{i+j-2q}\left(\frac{1}{i}\binom{i}{2q} + \frac{1}{j}\binom{j}{2q}\right) \right. \\
 + \binom{i+k-2q}{2t-2q}\frac{1}{i+k-2q}\left(\frac{1}{i}\binom{i}{2q} + \frac{1}{k}\binom{k}{2q}\right) \\
 \left. + \binom{j+k-2q}{2t-2q}\frac{1}{j+k-2q}\left(\frac{1}{j}\binom{j}{2q}+\frac{1}{k}\binom{k}{2q}\right)
 \right]\mathcal{B}_{i+j+k-2t}(x)  \\
 {-} \frac{(-1)^j}{ij\binom{i+j}{i}}B_{i+j}\mathcal{B}_k(x)
 {-} \frac{(-1)^k}{ik\binom{i+k}{i}}B_{i+k}\mathcal{B}_j(x) 
 {-} \frac{(-1)^k}{jk\binom{j+k}{j}}B_{j+k}\mathcal{B}_i(x) 
 - \frac{1}{2}\mathcal{B}_{i+j+k-2}(x) {+} \mathrm{const}.
\end{multline}

\begin{proof}
Let $i = n+1\geq2$, $j = m+1\geq2$, $k = l+1\geq2$ and $x = a+1$.
We calculate the LHS of~\eqref{RBnml} as 
\begin{multline}\label{RBnmlLHS}
(\mathcal{B}_{i}(x)-\mathcal{B}_{i}) 
(\mathcal{B}_{j}(x)-\mathcal{B}_{j})
(\mathcal{B}_{k}(x)-\mathcal{B}_{k}) 
 = \mathcal{B}_{i}(x)\mathcal{B}_{j}(x)\mathcal{B}_{k}(x) \\
 - (\mathcal{B}_{i}\mathcal{B}_{j}(x)\mathcal{B}_{k}(x)
  + \mathcal{B}_{j}\mathcal{B}_{i}(x)\mathcal{B}_{k}(x)
  + \mathcal{B}_{k}\mathcal{B}_{i}(x)\mathcal{B}_{j}(x)) \\
  + (\mathcal{B}_{i}\mathcal{B}_{j}\mathcal{B}_{k}(x)
  + \mathcal{B}_{i}\mathcal{B}_{k}\mathcal{B}_{j}(x)
  + \mathcal{B}_{j}\mathcal{B}_{k}\mathcal{B}_{i}(x)) 
  - \mathcal{B}_{i}\mathcal{B}_{j}\mathcal{B}_{k}.
\end{multline}

The last term on the RHS of~\eqref{RBnml} equals
\begin{equation}\label{RBnmlRHS3}
\mathcal{B}_{i+j+k-2}(x) - \mathcal{B}_{i+j+k-2}.
\end{equation}

We also have 
\begin{multline}\label{RBnmlRHS2}
- R(R(a^n)a^{m+l} + R(a^m)a^{n+l} + R(a^l)a^{n+m}) \\
 = -\sum\limits_{q\geq0}\left(\frac{1}{i}\binom{i}{2q}+\frac{1}{j}
 \binom{j}{2q}+\frac{1}{k}\binom{k}{2q}\right) B_{2q}\mathcal{B}_{i+j+k-1-2q}(x)  \\
 + \mathcal{B}_i\mathcal{B}_{j+k-1}(x) + \mathcal{B}_j\mathcal{B}_{i+k-1}(x)
 + \mathcal{B}_k\mathcal{B}_{i+j-1}(x) - \frac{3}{2}\mathcal{B}_{i+j+k-2}(x) + \mathrm{const}.
\end{multline}

We write down 
\begin{multline}\label{RBnmlRHS1}
R(R(R(a^n)a^m)a^l) 
 = \frac{1}{n+1}\sum\limits_{p=0}^n(-1)^{n-p}\binom{n+1}{n-p}B_{n-p}\frac{1}{m+p+2} \\
 \times\sum\limits_{s=0}^{p+m+1}(-1)^{m+p+1-s}\binom{m+p+2}{m+p+1-s}B_{m+p+1-s} 
 (\mathcal{B}_{l+s+2}(x)-\mathcal{B}_{l+s+2}). 
\end{multline}
We want to transform~\eqref{RBnmlRHS1} to the form
\begin{equation}\label{RBnmlRHS1-Trans}
\sum\limits_{q\geq0}\frac{1}{i}\binom{i}{2q}B_{2q}
\sum\limits_{t\geq0}\frac{1}{i+j-2q}\binom{i+j-2q}{2t}B_{2t}\mathcal{B}_{i+j+k-2t-2q}(x) + \mathrm{const}. 
\end{equation}

By exchange $n-p = 2q$ in~\eqref{RBnmlRHS1}, we get the summand
\begin{equation}\label{RBnmlRHS1-1}
\frac{1}{2}R(R(a^{n+m})a^l).
\end{equation}

Further, by exchange $m + p + 1 - s = 2t$, we have the additional summands
\begin{equation}\label{RBnmlRHS1-3}
\frac{1}{2}\sum\limits_{q\geq0}\frac{1}{i}\binom{i}{2q}B_{2q}(\mathcal{B}_{i+j+k-1-2q}(x)-\mathcal{B}_{i+j+k-1-2q})
 - \frac{1}{2}\mathcal{B}_i(\mathcal{B}_{j+k-1}(x)-\mathcal{B}_{j+k-1}).
\end{equation}

If we let $2q$ be equal~$i$, up constant term we get the summand 
\begin{equation}\label{RBnmlRHS1-2}
-\mathcal{B}_i\sum\limits_{t=0}^{[(j-1)/2]}\frac{1}{j}\binom{j}{2t}B_{2t}\mathcal{B}_{j+k-2t}(x)
 = -\mathcal{B}_i\sum\limits_{t\geq0}\frac{1}{j}\binom{j}{2t}B_{2t}\mathcal{B}_{j+k-2t}(x)
 + \mathcal{B}_i\mathcal{B}_j\mathcal{B}_k(x).
\end{equation}

Finally, letting $2t$ be equal~$j$, we get~\eqref{RBnmlRHS1-Trans} and the additional summand
\begin{equation}\label{RBnmlRHS1-4}
(\mathcal{B}_k-\mathcal{B}_k(x))\sum\limits_{q\geq0}\frac{1}{i}\binom{i}{2q}B_{2q}\mathcal{B}_{i+j-2q}.
\end{equation}

Applying the formula
$$
R(R(a^n)R(a^m)a^l)
 = R(R(R(a^n)a^m)a^l)
 + R(R(R(a^m)a^n)a^l)
 - R(R(a^{n+m})a^l),
$$
summing all such six expressions, by the formulas~\eqref{RBnmlLHS}--\eqref{RBnmlRHS1-4} we prove the statement. 
We have rewritten the sum of~\eqref{RBnmlRHS1-4} and the analogue of~\eqref{RBnmlRHS1-4} for~$j$ by~\eqref{Nielsen}, 
the sum equals
$$
(\mathcal{B}_k-\mathcal{B}_k(x))\left( \mathcal{B}_i\mathcal{B}_j + \frac{(-1)^i (i-1)!(j-1)!}{(i+j)!}B_{i+j}\right).
$$
Theorem is proved. \end{proof}

{\bf Corollary 1}. 
For all $i,j,k>0$, we have
\begin{multline}\label{SingleTripleBernIntegr}
\int\limits_{0}^x\mathcal{B}_i(y)\mathcal{B}_j(y)\mathcal{B}_k(y)\,dy 
= \sum\limits_{q,t\geq0}B_{2q}B_{2t-2q} \left[   
   \binom{i+j-2q}{2t-2q}\frac{1}{i+j-2q}\left(\frac{1}{i}\binom{i}{2q} + \frac{1}{j}\binom{j}{2q}\right) \right. \\
 + \binom{i+k-2q}{2t-2q}\frac{1}{i+k-2q}\left(\frac{1}{i}\binom{i}{2q} + \frac{1}{k}\binom{k}{2q}\right) \\
 \left. + \binom{j+k-2q}{2t-2q}\frac{1}{j+k-2q}\left(\frac{1}{j}\binom{j}{2q}+\frac{1}{k}\binom{k}{2q}\right)
 \right]\frac{\mathcal{B}_{i+j+k+1-2t}(x)}{i+j+k-2t}  \\
 - \frac{1}{ijk}\left(\frac{(-1)^j}{\binom{i+j}{i}}B_{i+j}\mathcal{B}_{k+1}(x)
 + \frac{(-1)^k}{\binom{i+k}{i}}B_{i+k}\mathcal{B}_{j+1}(x) 
 + \frac{(-1)^k}{\binom{j+k}{j}}B_{j+k}\mathcal{B}_{i+1}(x)\right) \\
 - \frac{\mathcal{B}_{i+j+k-1}(x) }{2(i+j+k-2)} + \mathrm{const}. 
\end{multline}

\begin{proof}
It follows from~\eqref{SingleTripleBern} and the formula
$\int\limits_{0}^xB_n(y)\,dy = \mathcal{B}_{n+1}(x)-\mathcal{B}_{n+1}$.
\end{proof}

Let us denote by $H_{n,s} = 1+1/2^s+\ldots+1/n^s$. So, $H_n = H_{n,1}$.

{\bf Corollary 2}. 
For all $N\geq3$, we have 
\begin{multline}\label{TripleSumDiv}
\frac{1}{3!}\sum\limits_{i+j+k = N;i,j,k>0}\mathcal{B}_i(x)\mathcal{B}_j(x)\mathcal{B}_k(x) \\ 
 = \sum\limits_{t>0}\binom{N-1}{2t}\mathcal{B}_{N-2t}(x)\left(
 \mathcal{B}_{2t}(H_{N-1}-H_{2t})
 + \frac{1}{2!}\sum\limits_{i+j=2t;i,j>0}\mathcal{B}_i\mathcal{B}_j  \right) \\
 - \frac{1}{12}\binom{N-1}{2}\mathcal{B}_{N-2}(x) + \frac{H_n^2 - H_{n,2}}{2}\mathcal{B}_N(x) + \mathrm{const}.
\end{multline}

\begin{proof}
We apply the formula~\eqref{SingleTripleBern} for all $i+j+k+N\geq3$ and $i,j,k>0$.
Let us compute the summation of the RHS of~\eqref{SingleTripleBern} divided by 6 for $q > 0$:
\begin{multline}\label{SumOfRBijk1-Main}
A = \sum\limits_{t\geq0,q>0}B_{2q}B_{2t-2q}
 \sum\limits_{i=1}^{N-2}\frac{1}{i}\binom{i}{2q}\sum\limits_{j=1}^{N-i-1}\frac{1}{i+j-2q}\binom{i+j-2q}{2t-2q}\mathcal{B}_{N-2t}(x) \\
 \allowdisplaybreaks
 = \sum\limits_{t>q>0}B_{2q}\mathcal{B}_{2t-2q}
 \sum\limits_{i=1}^{N-2}\frac{1}{i}\binom{i}{2q}\left(\binom{N-2q-1}{2t-2q}-\binom{i-2q}{2t-2q} \right) \mathcal{B}_{N-2t}(x) \\
 + \sum\limits_{t>0}B_{2t}\sum\limits_{i=1}^{N-2}\frac{1}{i}\binom{i}{2t}(H_{N-1-2q}-H_{i-2t})\mathcal{B}_{N-2t}(x) \\
 = \sum\limits_{t>q>0}\mathcal{B}_{2q}\mathcal{B}_{2t-2q}\mathcal{B}_{N-2t}(x)
 \left( \binom{N-1-2q}{2t-2q}\binom{N-2}{2q} - \binom{2t-1}{2q-1}\binom{N-2}{2t} \right) \\ 
  + \sum\limits_{t>0}\mathcal{B}_{2t}\left(\binom{N-2}{2t}H_{N-1-2t}-\binom{N-2}{2t}H_{N-2t-2}+\frac{1}{2t}\binom{N-2}{2t}-\frac{1}{2t}\right)\mathcal{B}_{N-2t}(x) \\
  = \sum\limits_{t>q>0}\binom{N-1}{2t}\binom{2t-1}{2q}\mathcal{B}_{2q}\mathcal{B}_{2t-2q}\mathcal{B}_{N-2t}(x) \\ 
 + \sum\limits_{t>0}\frac{1}{2t}\binom{N-1}{2t}\mathcal{B}_{2t}\mathcal{B}_{N-2t}(x)
 - \sum\limits_{t>0}\frac{1}{2t}\mathcal{B}_{2t}\mathcal{B}_{N-2t}(x) \\
  = \sum\limits_{t>0}\binom{N-1}{2t}\mathcal{B}_{N-2t}(x)\left(
 \frac{\mathcal{B}_{2t}}{2t} + \sum\limits_{t>q>0}\binom{2t-1}{2q}\mathcal{B}_{2q}\mathcal{B}_{2t-2q} \right)
 - \sum\limits_{t>0}\frac{1}{2t}\mathcal{B}_{2t}\mathcal{B}_{N-2t}(x) \\
 = \sum\limits_{t>0}\binom{N-1}{2t}\mathcal{B}_{N-2t}(x) 
  \sum\limits_{q>0}\frac{1}{2t}\binom{2t}{2q}\mathcal{B}_{2q}B_{2t-2q} 
 - \sum\limits_{t>0}\frac{1}{2t}\mathcal{B}_{2t}\mathcal{B}_{N-2t}(x).
\end{multline}
By the Miki identity~\eqref{Miki}, we may rewrite~\eqref{SumOfRBijk1-Main} as 
\begin{equation}\label{SumOfRBijk1-Main-2}
A = \sum\limits_{t>0}\binom{N-1}{2t}\mathcal{B}_{N-2t}(x)\left(-H_{2t-1}\mathcal{B}_{2t} + \frac{1}{2!}\sum\limits_{\substack{i+j=2t, \\ i,j>0}}\mathcal{B}_i\mathcal{B}_j\right)
 - \sum\limits_{q>0}\frac{1}{2q}\mathcal{B}_{2q}\mathcal{B}_{N-2q}(x).
\end{equation}
Above we have used the equalities
\begin{multline*}
\sum\limits_{j=1}^{N-i-1}\frac{1}{i+j-2q}\binom{i+j-2q}{2t-2q}
 = \frac{1}{2t-2q}\sum\limits_{j=1}^{N-i-1}\binom{i+j-2q-1}{2t-2q-1} \\
 = \frac{1}{2t-2q}\sum\limits_{s=i-2q}^{N-2q-2}\binom{s}{2t-2q-1}
 = \frac{1}{2t-2q}\left(\binom{N-2q-1}{2t-2q}-\binom{i-2q}{2t-2q} \right);
\end{multline*}
\begin{multline*}
\sum\limits_{i=1}^{N-2}\binom{i-1}{2q-1}H_{i-2q}
 = \sum\limits_{j=1}^{N-2-2q}\frac{1}{j}\sum\limits_{i=j+2q}^{N-2}\binom{i-1}{2q-1} \\ \allowdisplaybreaks
 = \sum\limits_{j=1}^{N-2-2q}\frac{1}{j}\left(\binom{N-2}{2q}-\binom{j+2q-1}{2q}\right) \\
 = \binom{N-2}{2q}H_{N-2q-2}-\frac{1}{2q}\sum\limits_{j=1}^{N-2-2q}\binom{j+2q-1}{2q-1} \\ 
 = \binom{N-2}{2q}H_{N-2q-2}-\frac{1}{2q}\left(\binom{N-2}{2q}-1\right). 
\end{multline*}

The summation of the RHS of~\eqref{SingleTripleBern} for $q = 0$ gives us
\begin{multline}\label{SumOfRBijk1-Add}
3\sum\limits_{t\geq0}B_{2t}
 \sum\limits_{k=1}^{N-2}\sum\limits_{i+j=N-k;i,j>0}\binom{i+j}{2t}\frac{1}{i+j}\left(\frac{1}{i}+\frac{1}{j}\right)\mathcal{B}_{N-2t}(x) \\
 = 3\sum\limits_{t\geq0}B_{2t}\sum\limits_{k=1}^{N-2}\binom{N-k}{2t}
 \sum\limits_{i=1}^{N-k-1}\frac{1}{N-k}\left(\frac{1}{i}+\frac{1}{N-k-i}\right) \mathcal{B}_{N-2t}(x) \\ 
 \allowdisplaybreaks
 = 6\sum\limits_{t>0}B_{2t}\sum\limits_{k=1}^{N-2}\binom{N-k}{2t}\frac{H_{N-k-1}}{N-k}\mathcal{B}_{N-2t}(x) 
 + 6\mathcal{B}_N(x)\sum\limits_{k=1}^{N-2}\frac{H_{N-k-1}}{N-k}\\
 = 6\sum\limits_{t>0}\mathcal{B}_{2t}\binom{N-1}{2t}\left(H_{N-1}-\frac{1}{2t}\right)\mathcal{B}_{N-2t}(x)
 + 3\big(H_{N-1}^2-H_{N-1,2}\big)\mathcal{B}_N(x).
\end{multline}
Here we have used the formulas (see~\cite[p.~279--280]{GKP})
$$
\sum\limits_{s=1}^{n-1}\binom{s}{m}H_s = \binom{n}{m+1}\left(H_n-\frac{1}{m+1}\right),\quad
\sum\limits_{s=1}^n\frac{H_s}{s} = \frac{1}{2}\big(H_n^2+H_{n,2}\big).
$$

The sum of first three summands from the fourth line of~\eqref{SingleTripleBern} by~\eqref{AlternOverBinom} gives
\begin{multline}\label{SumOfRBijk2}
-3\sum\limits_{k=1}^{N-2}B_{N-k}\mathcal{B}_k(x)\frac{1}{(N-k)(N-k-1)}\sum\limits_{i=1}^{N-k-1}\frac{(-1)^i}{\binom{N-k-2}{i-1}} \\
 = 3\sum\limits_{k=1}^{N-2}\frac{(1+(-1)^{N-k})}{(N-k)(N-k-1)}\frac{N-k-1}{N-k}B_{N-k}\mathcal{B}_k(x) 
 = 6\sum\limits_{0<t<N/2}\frac{1}{2t}\mathcal{B}_{2t}\mathcal{B}_{N-2t}(x).
\end{multline}

Finally, for the last nonconstant term of~\eqref{SingleTripleBern} we have
\begin{equation}\label{SumOfRBijk3}
\sum\limits_{i=1}^{N-2}(N-1-i)\mathcal{B}_{N-2}(x)
 = \binom{N-1}{2}\mathcal{B}_{N-2}(x).
\end{equation}

The formulas~\eqref{SumOfRBijk1-Main-2}--\eqref{SumOfRBijk3} imply~\eqref{TripleSumDiv}.
\end{proof}

In~\cite{KimProducts}, authors studied the expressions
\begin{equation}\label{SumDef}
S_{\geq0}(r,n) = \sum\limits_{i_1+\ldots+i_r = n}B_{i_1}(x)\ldots B_{i_r}(x),
\end{equation}
where the sum runs over all nonnegative integers $i_1,\ldots,i_r$ and $r\geq1$.
Moreover, in~\cite{KimProducts} it was shown that
\begin{equation}\label{SumDecomp}
S_{\geq0}(r,n)
 = \sum\limits_{k=1}^n B_k(x)\left(\frac{\binom{n+r}{k}}{n+r} 
 \sum\limits_{i_1+\ldots+i_r = n-k+1}B_{i_1}(1)\ldots B_{i_r}(1)-B_{i_1}\ldots B_{i_r}\right) + \mathrm{const},
\end{equation}
where there is a formula for the last constant term.
The simple idea lies behind such kind of decomposition:
polynomials $1 = B_0(x), B_1(x),\ldots,B_n(x)$ form a linear basis
of the vector space $\mathrm{Span}(1,x,x^2,\ldots,x^n)$.
Also, $(B_k(x))' = kB_{k-1}(x)$, the property which is not fulfilled 
for the polynomials $1,\mathcal{B}_1(x),\ldots,\mathcal{B}_n(x)$,
and so the approach from~\cite{KimProducts} could not be applied directly
for calculating the analogues of $S_{\geq0}(r,n)$ or $S_{>0}(r,n)$
for the sum of products $\mathcal{B}_{i_1}(x)\ldots\mathcal{B}_{i_r}(x)$.

Let us decompose $S_{\geq0}(r,n) = \sum\limits_{k=0}^n \alpha_k B_{n-k}(x)$.
In~\cite{KimProducts}, it was shown that $\alpha_{n-1} = 0$.
It is easy to see from~\eqref{SumDecomp} that all odd coefficients $\alpha_{2s+1}$ are zero.

{\bf Remark 4}.
The formula~\eqref{SumDecomp} applied for $r = 2$ up to constant coincides with~\eqref{Agoh14}.

Due to~\eqref{SumDecomp}, we have
\begin{multline}\label{Kim3}
S_{\geq0}(3,n) 
 = \sum\limits_{k=1}^n B_k(x)\left(\frac{\binom{n+3}{k}}{n+3} 
 \sum\limits_{a+b+c = n-k+1}B_{a}(1)B_b(1)B_c(1)-B_a B_b B_c\right) \\
 = -\frac{2}{n+3}\sum\limits_{t\geq0}\binom{n+3}{n-2t}\left(\sum\limits_{a+b+c=2t+1}B_a B_b B_c\right)B_{n-2t}(x) \\
 \!\!{=} \binom{n+2}{2}B_n(x) {+} \frac{1}{4}\binom{n+2}{4}B_{n-2}(x)
 + \frac{3}{n+3}\sum\limits_{t\geq2}\binom{n+3}{n-2t}B_{n-2t}(x)\sum\limits_{q=0}^{t}B_{2q}B_{2t-2q}.
\end{multline}
If we calculate the sum $S_{\geq0}(3,n)$ due to the formula~\eqref{SingleTripleBern}, 
we will get the same as~\eqref{Kim3} modulo the equality~\eqref{Agoh14}.

\section*{Acknowledgments}
Author is supported by the Austrian Science Foun\-da\-tion FWF, grant P28079.

Author is grateful to Oleg Ogievetsky and Li Guo for pointing out 
the references~\cite{Ogievetsky} and \cite{Miller} respectively.

\medskip

\noindent Vsevolod Gubarev \\
University of Vienna \\
Oskar-Morgenstern-Platz 1, 1090, Vienna, Austria \\
e-mail: vsevolod.gubarev@univie.ac.at

\end{document}